\documentclass[10pt]{article}
\usepackage{hyperref}
\usepackage{fancyvrb}  
\usepackage{amsmath,amssymb,color,bbm,ifthen,enumerate}
\usepackage{graphicx}

\newtheorem{thm}{Theorem}[section]

\newtheorem{prop}[thm]{Proposition}
\newtheorem{cor}[thm]{Corollary}

\newenvironment{rem}[1][Remark.]{\begin{trivlist}
\item[\hskip \labelsep {\bfseries #1}]}{\end{trivlist}}

\newenvironment{pf}{\paragraph{Proof.}}
{\nopagebreak\hfill\nopagebreak\rule{2mm}{2mm}\par\bigskip}

\newcommand{\bone}{\mbox{\boldmath$1$}}

\newcommand{\argst}{\{ s_i,t_i \}_{i=1}^n }
\newcommand{\argsl}{\{ s_i,\lambda_i \}_{i=1}^n }

\begin{document}

\renewcommand{\baselinestretch}{1.3}

\title{Using Differential Equations to Obtain Joint Moments of
First-Passage Times of Increasing L\'{e}vy Processes\footnote{Revised} }

\author{ Mark Veillette  and Murad S. Taqqu\thanks{ This research was partially supported by the NSF grants DMS-0505747,
DMS-0706786, and DGE-0221680.}
\\ Boston University}

\maketitle

\begin{abstract}
Let $\{D(s), s \geq 0 \}$ be a L\'{e}vy subordinator, that is, a
non-decreasing process with stationary and independent increments and
suppose that $D(0) = 0$.  We study the first-hitting time of the
process $D$, namely, the process $E(t) = \inf \{ s: D(s) > t \}$, $t
\geq 0$.  
 The process $E$ is, in general,  non-Markovian with non-stationary
 and non-independent
increments.   We derive a partial
differential equation for the Laplace transform of the $n$-time tail
distribution function $P[E(t_1) > s_1,\dots,E(t_n) > s_n]$, and show
that this PDE has a unique solution given natural boundary conditions.
This PDE can be used to derive all $n$-time moments
of the process $E$.   
\end{abstract}

\section{Introduction}

Consider a non-decreasing L\'{e}vy Process $\{D(s),\ s \geq 0 \}$,
starting from $0$, which is continuous from the right with left
limits.  Such a process is called a subordinator.  It has stationary
and independent increments and is characterized by its Laplace Transform
\begin{equation}
 \mathbb{E} e^{-\lambda D(s)} = e^{-s \phi(\lambda)}, \quad \lambda
 \geq 0. \nonumber
\end{equation}
The function $\phi$ is called the Laplace exponent and is given by the L\'{e}vy-Khintchine formula:
\begin{equation}\label{e:LKformula}
\phi(\lambda) = \mu \lambda + \int_{(0,\infty)} \left( 1 - e^{-\lambda x} \right) \Pi(dx),
\end{equation}
where $\mu \geq 0$ is the drift and $\Pi$ is a measure on
$\mathbb{R}^+ \cup \{0 \}$ which satisfies $\int_0^\infty (1 \wedge x)
\Pi(dx) < \infty$ (see ~\cite{Applebaum:2004},
~\cite{Bertoin:1996} or ~\cite{Bertoin:1999}). 

We want to study the first-passage time of such a
process, focusing on its finite-dimensional distributions.  The
first-passage time of a subordinator $\{D(s), s \geq 0\}$, is a new
process $\{E(t), t \geq 0\}$, commonly called an \textit{inverse
  subordinator}, and is defined as follows:
\begin{equation}
E(t) = \inf\{ s: D(s) > t \},\quad t \geq 0. \nonumber
\end{equation}   
 It is worth noting that in certain
cases, the L\'{e}vy subordinator $D$ is itself a first passage time of
another L\'{e}vy processes.  For example, the $\frac{1}{2}$-stable
subordinator is the first
passage time of standard Brownian motion, and the inverse Gaussian
subordinator is the first passage time of standard Brownian
motion with drift (\cite{Applebaum:2004}, exercise
2.2.10).  For these examples, $E$ is the
``first-passage time of a first-passage time''.

Inverse subordinators appear in a variety of applications.  For
instance, they are used
in the study of 
scaling limits of continuous-time random walks and fractional kinetics,
~\cite{Baule:2005},~\cite{Baule:2007},~\cite{Meerschaert:2006},~\cite{Meerschaert:2004},
  ~\cite{Arous:2007}. Here, the time variable of a Markov process (typically Brownian
  motion), is replaced by an inverse $\alpha$-stable
  subordinator, with $0 <\alpha < 1$. This particular time change gives rise to anomalous diffusion,
  or sub-diffusion,
  where the variance of the process grows at a rate which is
  non-linear in time.  The $\alpha$-stable subordinator is ``fast'' because
  all moments of order $\alpha$ or less are infinite.  An ``ultrafast''
  subordinator introduced in \cite{Meerschaert1:2004} and \cite{Meerschaert:2006} is a
  process which has, in general, no finite moments.  While these
  subordinators have infinite moments, the moments of their inverses
  are finite.  For the $\alpha$-stable subordinator, the first moment of its
  inverse subordinator, $\mathbb{E}E(t)$, grows as $t^\alpha$.  For the
  ultrafast subordinator, the inverse is ``ultraslow'', and can grow like
  a slowly varying function.

Our goal in this paper is to characterize the $n$-point distribution
function of a general inverse subordinator with a simple partial
differential equation and to
use this PDE to derive explicit expressions for joint moments.  A PDE
was derived (heuristically) by ~\cite{Baule:2005},~\cite{Baule:2007}
for the joint density function in the special case of $\alpha$-stable
inverse subordinators and was used to obtain joint moments of these
processes.  To obtain precise results in the general case, we found it
convenient to use suitable modified cumulative distribution functions
instead of densities.  An alternative approach, using Cox processes
was followed by  ~\cite{Nord:2005}, who also obtains joint moments of
increments of inverse subordinators.  Our results are implemented in
\cite{veillette1:2008}, where we develop algorithms for computing the
moments numerically.

 This paper is organized as follows:   The $n$-point
distribution of $E$ is studied in Section \ref{s:multipletimes},  where we
derive a PDE for the Laplace transform of the $n$-point tail
distribution function $P[E(t_1)>s_1,\dots,E(t_n) > s_n]$.  In Section
\ref{s:moments}, we use this PDE to calculate the $n$-point
moments of a general inverse subordinator. The examples presented briefly in Section
\ref{s:exbrief} are discussed in detail in the companion paper
\cite{veillette1:2008}.  For convenience,  we provide an appendix
which contains key path properties of inverse subordinators.

\section{PDE for the Multiple-time Distribution Function of Inverse Subordinators}\label{s:multipletimes}

Unlike the L\'{e}vy process $D$, the inverse subordinator $E$ is often non-Markovian.  Thus, we must consider {\it all} finite-dimensional distributions.  For
$n \geq 1$, define
\begin{equation}\label{e:tail}
H^{(n)}(s_1,\dots,s_n,t_1,\dots,t_n) = P[E(t_i) > s_i,\
i=1,\dots,n].
\end{equation}  
  To simplify notation, we will write
$H^{(n)}(\{s_i,t_i \}_{i=1}^n) =
H^{(n)}(s_1,\dots,s_n,t_1,\dots,t_n)$.  In the following,
$\mathbb{R}^n_+ = \{(s_1,\dots,s_n), \ s_i \geq 0, i=1,\dots,n \}$.
Also, $C^1(\mathbb{R}^n_+)$ will denote the space of continuously
differentiable functions on $\mathbb{R}^n_+$.  The following theorem
shows the Laplace transform, 
\begin{equation}
\widetilde{H^{(n)}}(\argsl) = \int_0^\infty \dots \int_0^\infty
e^{-\lambda_1 t_1 - \dots - \lambda_n t_n} H^{(n)}(\argst)
dt_1\dots dt_, \nonumber
\end{equation}
of $H^{(n)}$ satisfies a PDE.  Observe that
$\widetilde{H^{(1)}}(0,\lambda) = \int_0^\infty P[E(t) > 0] e^{-\lambda
  t} dt = 1/\lambda$ since $P[E(t) > 0] = 1$ for $t>0$.  We let $H^{(n)}(\{s_i = 0 \})$ denote
$H^{(n)}$ in (\ref{e:tail}) with $s_i=0$ and the other arguments
unchanged.

\begin{thm}\label{t:taildist}
Let $D$ be a general L\'{e}vy subordinator and let $E$ be the inverse
subordinator of $D$. For $s_1, s_2, \dots,
s_n \geq 0$, the Laplace Transform of the $n$-point tail distribution $H^{(n)}$ of $E$ defined by
(\ref{e:tail}) is the unique solution in $C^1(\mathbb{R}^n_+)$ to the following PDE
\begin{equation}\label{e:inversePDE}
\left( \frac{\partial}{\partial s_1} + \frac{\partial}{\partial s_2} +
  \dots + \frac{\partial}{\partial s_n} \right) \widetilde{H^{(n)}} = -\phi(\lambda_1 +
\dots + \lambda_n) \widetilde{H^{(n)}},
\end{equation}
together with the boundary conditions 

\begin{equation}\label{e:inverseBCs}
\begin{cases}
\widetilde{H^{(1)}}(0,\lambda) = \frac{1}{\lambda}, \quad n=1, \\
\widetilde{H^{(n)}}(\{s_i=0\}) =
\frac{1}{\lambda_i}\widetilde{H^{(n-1)}}(\{s_k,\lambda_k\}_{k\neq i}),
\quad i=1,\dots,n, \quad n>1. \end{cases}
\end{equation}
\end{thm}


\begin{pf}

From Proposition \ref{p:setequal}, we have that 
\begin{eqnarray*}
H^{(n)}(\argst) &=&  P[E(t_i) > s_i,\
i=1,\dots,n]. \\
&=&  P[D(s_i) < t_i,\
i=1,\dots,n],\quad \mbox{for a.e.} \ t_1,\dots,t_n \geq 0.
\end{eqnarray*}
Let $ s_0=0$ and fix $ s_1,  s_2,\dots, s_n \geq 0$.  Calculating the Laplace transform
in $t_1,\dots,t_n$ of $H^{(n)}$, we have
\begin{equation}
\widetilde{H^{(n)}}(\argsl) =
\int_0^\infty \dots \int_0^\infty \exp(-\sum_{i=1}^n \lambda_i t_i)
H^{(n)}(\argst) dt_1\dots dt_n \nonumber
\end{equation}
\begin{equation}
= \frac{(-1)^n}{\prod_{i=1}^n \lambda_i } \int_0^\infty \dots
\int_0^\infty \left( \frac{\partial^n}{\partial t_1 \dots \partial t_n}
\exp(-\sum_{i=1}^n \lambda_i t_i) \right)  P[D(s_i) < t_i,\
i=1,\dots,n]
 dt_1\dots dt_n. \nonumber
\end{equation}
Integrating by parts in $t_1,\dots,t_n$, we obtain
\begin{equation}
\widetilde{H^{(n)}}(\argsl)  =
\frac{(-1)^{2n}}{\prod_{i=1}^n \lambda_i } \mathbb{E}
\exp(-\sum_{i=1}^n \lambda_i D(s_i)). \nonumber
\end{equation}

Now rewrite the variables $s_0,s_1,\dots,s_n$ is an increasing order,
$0 = s_{j(0)} \leq s_{j(1)} \leq \dots \leq s_{j(n)}$, where $j(0)=0$ and
$j(1),\dots,j(n)$ is a
permutation of the integers $1,\dots,n$.   
We add and subtract terms in the
summation above and use the independence and stationarity of the
increments of $D$
to obtain
\begin{eqnarray}
\widetilde{H^{(n)}}(\argsl)  &=&
\frac{1}{\prod_{i=1}^n \lambda_i } \mathbb{E} \exp \left( -\sum_{i=1}^n
\left(\sum_{k=i}^n \lambda_{j(k)} \right) (D(s_{j(i)}) - D(s_{j(i-1)})) \right)\nonumber \\
&=&  \frac{1}{\prod_{i=1}^n \lambda_i } \prod_{i=1}^n \mathbb{E} \exp
\left( -
\left(\sum_{k=i}^n \lambda_{j(k)} \right) D(s_{j(i)} - s_{j(i-1)})
\right) \nonumber \\
&=& \frac{1}{\prod_{i=1}^n \lambda_i  } \exp \left( - \sum_{i=1}^n \phi
  \left( \sum_{k=i}^n \lambda_{j(k)} \right) (s_{j(i)} -
  s_{j(i-1)}) \right) \label{e:formofht},
\end{eqnarray}
Differentiating with respect to $s_i$,
$i=1,\dots,n$, we have that $\widetilde{H^{(n)}}$ satisfies the system
of PDEs 
\begin{eqnarray*}
\frac{\partial}{\partial s_{j(i)}} \widetilde{H^{(n)}} &=&
-\left(\phi( \sum_{k=i}^n \lambda_{j(k)}) - \phi(\sum_{k=i+1}^{n}
  \lambda_{j(k)}) \right) \widetilde{H^{(n)}}, \quad i=1,2,\dots,n-1, \\
\frac{\partial}{\partial s_{j(n)}} \widetilde{H^{(n)}} &=& -\phi(\lambda_{j(n)})
\widetilde{H^{(n)}}.
\end{eqnarray*}
Adding these $n$ equations together gives
\begin{equation}
\left( \frac{\partial}{\partial s_{j(1)}} + \dots +
  \frac{\partial}{\partial s_{j(n)}} \right)  \widetilde{H^{(n)}} =
-\phi( \lambda_{j(1)} + \dots + \lambda_{j(n)} ) \widetilde{H^{(n)}}. \nonumber
\end{equation}
Since $j$ is simply a permutation of the integers $\{ 1,\dots,n \}$,
this is equivalent to (\ref{e:inversePDE}).

Now we claim that with the boundary conditions (\ref{e:inverseBCs}), the
tail probability $\widetilde{H^{(n)}}$ is the unique solution to the
PDE (\ref{e:inversePDE}) in $C^1(\mathbb{R}^n_+)$.  Suppose $\widetilde{H}_1(\argsl)$ is
another solution in $C^1(\mathbb{R}^n_+)$ to (\ref{e:inversePDE}) with the boundary conditions
(\ref{e:inverseBCs}).  Define $V(\argsl) = (\widetilde{H^{(n)}} -
\widetilde{H}_1)(\argsl)$.  By linearity, $V$ is also a solution to
(\ref{e:inversePDE}) with boundary conditions 
\begin{equation}\label{e:newinverseBCs}
V(\{s_i=0\}) = 0,
\quad i=1,\dots,n.
\end{equation}
Fix $(s_1,\dots,s_n) \in \mathbb{R}^n_+$, and WLOG, assume $s_1 \leq
\dots \leq s_n$ (if not we can simply re-index as before).  Proceeding by the method of characteristics (\cite{Evans:1998},
section 3.2), we define the function 
\begin{equation}\label{e:defofnu_inverse}
\nu(\tau) = V(\{s_i +
\tau,\lambda_i \}_{i=1}^n) \quad \  \tau \geq -s_1.  
\end{equation}
Observe that (\ref{e:newinverseBCs}) can now be written as
\begin{equation}\label{e:rewriteBCs}
V(0,s_2,\dots,s_n,\lambda_1,\dots,\lambda_n) = 0.
\end{equation}
Differentiating (\ref{e:defofnu_inverse}) with respect
to $\tau$ and using the definition of $\nu$, we have
\begin{equation}\label{e:simpODE}
\nu'(\tau) = \left( \frac{\partial}{\partial s_1} + \frac{\partial}{\partial s_2} +
  \dots + \frac{\partial}{\partial s_n} \right) V(\{s_i +
\tau,\lambda_i\}_{i=1}^n). 
\end{equation} 
Since $V$ satisfies the PDE (\ref{e:inversePDE}), equation (\ref{e:simpODE}) implies that
$\nu$ satisfies the ODE
\begin{equation}
\nu'(\tau) = - \phi(\lambda_1+\dots+\lambda_n) \nu(\tau). \nonumber
\end{equation}
Now, using $\tau = -s_1$ as a initial condition and using
(\ref{e:rewriteBCs}), we have 
\begin{equation}
\nu(-s_1) =
V(0,s_2-s_1,\dots,s_n-s_1,\lambda_1,\dots,\lambda_n) = 0. \nonumber
\end{equation}
  Thus,
$\nu(\tau) \equiv 0$ for $\tau \geq -s_1$ is a solution to this initial
value problem, and by standard uniqueness theory of ODEs
(\cite{Cole:1968}, chapter 2), this is the
unique solution.  Hence 
\begin{equation}
0 = \nu(0) = V(\argsl) =  \widetilde{H^{(n)}}(\argsl) -
\widetilde{H}_1(\argsl), \nonumber
\end{equation}
 implying
$\widetilde{H^{(n)}}$ is the unique solution to (\ref{e:inversePDE})
with the boundary conditions (\ref{e:inverseBCs}).

\end{pf}

\begin{rem}
The boundary conditions (\ref{e:inverseBCs}) are natural ones, since
they imply that the
$n$-dimensional distribution function can be reduced to the $n-1$
dimensional distribution function when $s_i=0$ for some $i$. 
\end{rem}

\begin{rem}
One might wonder why a PDE for the $n$-time Laplace transform of $H^{(n)}$ is useful when in fact we can write down its solution (\ref{e:formofht}) in closed form.  It is, because:

\begin{itemize}

\item  It simplifies the calculation of the moments of $E$, which we demonstrate in Section \ref{s:moments}.

\item  Understanding the dynamics of $E$ is useful in the study of more complicated processes.  For example,  in \cite{Baule:2007}, the PDE (\ref{e:inversePDE}) in the case of the $\alpha$-stable process (where $\phi(\lambda) = \lambda ^\alpha$) is used to derive equations corresponding to the $n$-time distribution functions of the so-called {\it fractional kinetic} process $Z_{\alpha}(t) = B(E(t))$, where $B$ is Brownian motion and $E$ is an inverse $\alpha$-stable subordinator.  This process appears as a scaling limits for various trap models, \cite{Arous:2007}.  Thus, the PDE (\ref{e:inversePDE}) can extend this analysis done in \cite{Baule:2007} to a larger class of processes.    

\end{itemize}

\end{rem}

\section{Moments of Inverse Subordinators}\label{s:moments}

In this section we use Theorem \ref{t:taildist} to calculate moments
of a general inverse subordinator.  The utility of the PDE
(\ref{e:inversePDE}) will become clear in that it will simplify many
of the following computations.

Before we calculate moments, we first argue that all moments of an inverse subordinator are finite.  Notice that, for any $x > 0$, we can bound the tail distribution of $E$ using (\ref{e:inequality}) from the proof of Proposition \ref{p:setequal} and Markov's inequality:
\begin{equation}
P[E(t) > s]  \leq  P[D(s) \leq t] = P[ e^{-x D(s)} \geq e^{-x t}]  \leq e^{xt} \mathbb{E}e^{-x D(s)} = e^{x t} e^{-s \phi(x)},
\end{equation}
which implies that $\mathbb{E}E(t)^\gamma < \infty$ for any $\gamma > 0$.

 The Laplace transform of $\mathbb{E}E(t)^\gamma$ with $\gamma>0$ has a simple form:

\begin{eqnarray}
\mathbb{E}E(t)^\gamma  &=&  \gamma  \int_0^\infty e^{-\lambda t} \int_0^\infty  s^{\gamma-1} P[E(t) > s] ds  dt \nonumber\\
&=&  \gamma \int_0^\infty s^{\gamma-1} \int_0^\infty e^{-\lambda t}  (1- P[E(t)
\leq s]) dt ds \\
&=& \frac{\gamma}{\lambda} \int_0^\infty s^{\gamma-1} e^{-s \phi(s)}
ds \nonumber \\
&=& \frac{\gamma \Gamma(1+ \gamma)}{\lambda \phi(\lambda)^\gamma} \label{e:ltgammamoment},
\end{eqnarray}
where we have used integration by parts to calculate $\int_0^\infty e^{-\lambda t}  (1- P[E(t)
\leq s]) dt$.  Of particular importance is the mean of $E(t)$.  Let
$U(t) = \mathbb{E}E(t)$.  From (\ref{e:ltgammamoment}) $U$ has Laplace
transform given by
\begin{equation}\label{e:lt_mfpt}
\widetilde{U}(\lambda) = \frac{1}{\lambda \phi(\lambda)}.
\end{equation} 
We will see in the following that $U$ characterizes all
finite-dimensional distributions of the process $E$.  While the
Laplace transform of $U$ is easy to express in terms of $\phi$,
calculating the inverse is, in general, not always an easy task.  In
\cite{veillette1:2008}, a numerical method is given for computing $U$ for
general $\phi$.     

\begin{rem} Although the subordinator $D$ is, in general, different
  from a
  renewal process, it has some of its characteristics.  For example,
  it satisfies the  so called \textit{renewal theorem}, which states that if the
mean of the subordinator is finite, then $\mathbb{E}E(t) = U(t) \sim
\frac{t}{\mathbb{E}D(1)}$, as $t \rightarrow \infty$. This was proven
using various methods in ~\cite{Nord:2005}, ~\cite{Bertoin:1999} and ~\cite{Klafter:2003}.  This fact can be easily seen from the
Laplace transform of $U$ given above.  Indeed, if $D(1)$ has finite
mean, then from the L\'{e}vy-Khintchine formula, 
\begin{equation}\label{e:ratio}
\mathbb{E}D(1) = \phi'(0) = \lim_{\lambda \rightarrow 0}
\frac{\phi(\lambda)}{\lambda} =  \mu + \int_0^\infty x \Pi(dx) < \infty,
\end{equation}
and thus $\phi(\lambda) \sim \mathbb{E}D(1) \lambda$ as $\lambda
\rightarrow 0$.  From (\ref{e:lt_mfpt}), $\widetilde{U}(\lambda) \sim
1/(\mathbb{E}D(1) \lambda^2)$ as $\lambda \rightarrow 0$ and the Tauberian theorem
(\cite{Bertoin:1996}, page 10) implies the renewal theorem:
\begin{equation}
U(t) \sim \frac{t}{\mathbb{E}D(1)}, \quad t \rightarrow \infty.
\end{equation}
Thus, for subordinators with finite mean, their mean first-passage
time will exhibit a non-linear
transient behavior for small times, followed by a linear behavior for
large times.   
\end{rem}

As an application of the differential equation given in Theorem
\ref{t:taildist}, we obtain expressions of the Laplace transforms
for the $n$-time integer moments of $E$ in terms of Laplace transforms
of lower-order moments.
\begin{thm}\label{t:inversemoments}
Let $D$ be a general L\'{e}vy subordinator with L\'{e}vy exponent $\phi$ and let $E$ be the inverse
subordinator of $D$.  For positive integers $m_1,\dots,m_n$, let
\begin{equation}\label{e:momentnotation}
U(t_1,\dots,t_n;m_1,\dots,m_n) =
\mathbb{E}E(t_1)^{m_1}\dots E(t_n)^{m_n}.
\end{equation}
In the special case $n=1,m_1=1$, we will simply write $U(t,1) = U(t)$.  The $n$-time Laplace Transform
of $U$ is given in terms of strictly lower order moments by
\begin{equation}\label{e:inversenmoments}
\widetilde{U}(\lambda_1,\dots,\lambda_n;m_1,\dots,m_n) =
\frac{1}{\phi(\lambda_1+\dots+\lambda_n)} \sum_{i=1}^n m_i
\widetilde{U}(\lambda_1,\dots,\lambda_n;m_1,\dots,m_{i-1},m_i-1,m_{i+1},\dots,m_n).
\end{equation}
\end{thm}

\begin{pf}


We first write $U$ in terms of the tail probability $H^{(n)}$, 
\begin{equation}
U(t_1,\dots,t_n;m_1,\dots,m_n) = M  \int_0^\infty\dots \int_0^\infty
s_1^{m_1-1} \dots s_n^{m_n-1} H^{(n)}( \argst )
ds_1,\dots,ds_n, \nonumber
\end{equation}
where $M = \prod_{i=1}^n m_i$.  Taking the Laplace transform in $t_1,\dots,t_n$ and rearranging the order
of integration yields
\begin{equation}\label{e:momentstep1}
\widetilde{U}(\lambda_1,\dots,\lambda_n;m_1,\dots,m_n) = M \int_0^\infty\dots \int_0^\infty
s_1^{m_1-1} \dots s_n^{m_n-1} \widetilde{H^{(n)}}( \argsl )
ds_1,\dots,ds_n.
\end{equation}
Now multiply through by $\phi(\lambda_1+\dots + \lambda_n)$ and
apply Theorem \ref{t:taildist} to obtain
\begin{equation}
\phi(\lambda_1+\dots+\lambda_n) \widetilde{U}(\lambda_1,\dots,\lambda_n;m_1,\dots,m_n) \label{e:momentstep1.5}
\end{equation}
\begin{equation}
=-M\int_0^\infty\dots \int_0^\infty
s_1^{m_1-1} \dots s_n^{m_n-1}\left(   \frac{\partial}{\partial s_1} + 
  \dots + \frac{\partial}{\partial s_n} \right) \widetilde{H^{(n)}}( \argsl )
ds_1,\dots,ds_n. \nonumber
\end{equation}
\begin{equation}
=-M\sum_{i=1}^n \int_0^\infty \dots \int_0^\infty \left( \int_0^\infty
s_1^{m_1-1} \dots s_n^{m_n-1} \frac{\partial}{\partial s_i} \widetilde{H^{(n)}}( \argsl )
ds_i \right) ds_1 \dots ds_{i-1} ds_{i+1} \dots ds_n. \label{e:momentstep2}
\end{equation}

Let us now focus on the inner-most integral above.  If $m_i = 1$,
we have
\begin{eqnarray}
\int_0^\infty
s_1^{m_1-1} \dots s_n^{m_n-1} \frac{\partial}{\partial s_i} \widetilde{H^{(n)}}( \argsl )
ds_i &=& \prod_{k \neq i} s_k^{m_k-1} \int_0^\infty  \frac{\partial}{\partial s_i} \widetilde{H^{(n)}}( \argsl )
ds_i \nonumber \\
&=& - \prod_{k \neq i} s_k^{m_k-1} \frac{1}{\lambda_i}
\widetilde{H^{(n-1)}} ( \{s_k,\lambda_k \}_{k \neq i} ).\label{e:momentsm1}
\end{eqnarray}
Above we have used the fact that $\widetilde{H^{(n)}} \rightarrow 0$ exponentially
as $s_i \rightarrow \infty$ , which follows from (\ref{e:formofht}),
and that $\widetilde{H^{(n)}}(\{s_i=0\}) =\frac{1}{\lambda_i} \widetilde{H^{(n-1)}}(\{
s_k,\lambda_k\}_{k\neq i})$ from (\ref{e:inverseBCs}).  

If $m_i >1$, we integrate by parts and get
\begin{equation}
\int_0^\infty
s_1^{m_1-1} \dots s_n^{m_n-1} \frac{\partial}{\partial s_i} \widetilde{H^{(n)}}( \argsl )
ds_i = \prod_{k \neq i} s_k^{m_k-1} \int_0^\infty s_i^{m_i-1}  \frac{\partial}{\partial s_i} \widetilde{H^{(n)}}( \argsl )
ds_i \nonumber 
\end{equation}
\begin{equation}
= - (m_i-1) \prod_{k \neq i} s_k^{m_k-1} \int_0^\infty s_i^{m_i-2} \widetilde{H^{(n)}}( \argsl )
ds_i \label{e:momentsm2}
\end{equation}

With (\ref{e:momentsm1}) and (\ref{e:momentsm2}),
(\ref{e:momentstep1}) implies that (\ref{e:momentstep1.5}) and
(\ref{e:momentstep2}) can now be written as
\begin{equation}
\phi(\lambda_1+\dots+\lambda_n)
\widetilde{U}(\lambda_1,\dots,\lambda_n;m_1,\dots,m_n)  = 
\sum_{i=1}^n m_i
\widetilde{U}(\lambda_1,\dots,\lambda_n;m_1,\dots,m_{i-1},m_i-1,m_{i+1},\dots,m_n), \nonumber
\end{equation}
which finishes the proof.

\end{pf}

\begin{rem}
Theorem \ref{t:inversemoments} is equivalent to Theorem 2.1 in
\cite{Nord:2005}.  There, the result gives the joint
moments of increments of the inverse subordinator: $\mathbb{E} \prod_{i=1}^n (E(t_i) -
E(s_i))^{k_i}$, where $0 \leq s_1 \leq t_1 \leq \dots \leq s_n \leq
t_n$, and $k_i,\ i=1\dots n$ are positive integers, and expresses them
with an integral expression.  The formulas given
here, which are obtained using different methods, give directly
$\mathbb{E}E(t_1)^{m_1}\dots E(t_n)^{m_n}$ and have the added advantage of expressing higher order moments in
terms of lower
order ones.  Also, out results give the Laplace transform of all
$n$-time moments.    
\end{rem}

\medskip

\noindent \textbf{Recursion:} Theorem  \ref{t:inversemoments} gives expressions for the Laplace
transform of all $n$-time
moments of a general inverse subordinator. In theory, these
$n$-dimensional Laplace transforms can be inverted, however doing so directly can be a
formidable task.  As an alternative, we take advantage of the
recursive nature of equation (\ref{e:inversenmoments}).  Let $N \in
\mathbb{Z}$ be the order of the moment $U(t_1,\dots,t_n,m_1,\dots,m_n)$,
i.e. $N= m_1+m_2+\dots+m_n$.  Observe that equation
(\ref{e:inversenmoments}) implies that moments of order $N$
can be calculated using a linear combination of convolutions involving moments of order
$N-1$ and the inverse Laplace transform of the function $1/\phi$,
which we address below.
Thus, if one has $U(t) = U(t;1)$, then all moments can be obtained
inductively using this method.


\subsection{Getting the inverse Laplace transform}

To calculate moments of inverse subordinators, we must be able
to invert Laplace transforms of the form $\widetilde{f}/\phi$, where
the function $f$ is known. Hence, we must first
obtain the inverse Laplace transform of the function $1/\phi$.  Following \cite{Bertoin:1996},
section III.1, define the \textit{renewal measure} to be the Borel
measure whose distribution function is given by $U(t)$.  From this
we see that for a.e. $0 \leq a < b$, 
\begin{eqnarray}
 \int_0^\infty \bone_{(a,b]}(\tau) dU(\tau) &=& U(b) - U(a) \nonumber \\
&=& \int_0^\infty P[E(b) > s] - P[E(a) > s] ds \nonumber \\
&=& \int_0^\infty P[D(s) \leq b] - P[D(s) \leq a] ds \nonumber \\
&=& \mathbb{E} \int_0^\infty \left( \bone_{(-\infty,b]}(D(s)) -
  \bone_{(-\infty,a]}(D(s)) \right)
ds \nonumber \\
&=& \mathbb{E} \int_0^\infty \bone_{(a,b]}(D(s)) 
ds. \label{e:renewalmeasure}
\end{eqnarray} 
By approximating with step functions, this can be extended to
$\int_0^\infty g(\tau) dU(\tau) = \mathbb{E} \int_0^\infty g(D(s)) ds$
where $g$ is a continuous function.  Choosing $g(\tau) =
e^{-\lambda \tau}$, we get that the Laplace transform of the renewal
measure is given by
\begin{equation}
\int_0^\infty e^{-\lambda \tau} dU(\tau) = \mathbb{E} \int_0^\infty
e^{-\lambda D(s)} ds  
= \int_0^\infty e^{-s \phi(\lambda) }ds 
= \frac{1}{\phi(\lambda)}. \nonumber
\end{equation}

Thus, for an arbitrary function $f$, the inverse Laplace transform of
$\widetilde{f}/\phi$ is given by the convolution of $f$ with the
renewal measure, i.e. 
\begin{equation}\label{e:inversef1d}
{\cal{L}}^{-1}\left[\frac{1}{\phi(\lambda)}
  \widetilde{f}(\lambda)\right] (t) =
\int_0^t f(t-\tau) dU(\tau). 
\end{equation}
This can be generalized to $n$
dimensions. For ease of notation, write $dU(\tau) = U'(\tau)d\tau$,
where $U'$ is interpreted as a generalized function(since $U$ might
contain jumps).  If $0 \leq t_1 \leq t_2 \leq \dots
\leq t_n$, then
\begin{equation}
{\cal{L}}^{-1} \left[  \frac{1}{\phi(\lambda_1 + \lambda_2 + \dots +
  \lambda_n)} \widetilde{f}(\lambda_1,\lambda_2,\dots,\lambda_n)
\right] (t_1,\dots,t_n) \nonumber
\end{equation}
\begin{equation}
= \int_0^{t_1} \int_0^{t_2}\dots\int_0^{t_n}
f(t_1-\tau_1,\dots,t_n-\tau_n){\cal{L}}^{-1}\left[
  \frac{1}{\phi(\lambda_1+\dots+\lambda_n)}\right](\tau_1,\dots,\tau_n)
d\tau_1,\dots d\tau_n \nonumber
\end{equation}
\begin{equation}
= \int_0^{t_1} \int_0^{t_2} \dots \int_0^{t_n}
f(t_1-\tau_1,\dots,t_n-\tau_n) U'(\tau_1)\delta(\tau_2-\tau_1)\dots\delta(\tau_n-\tau_1)
d\tau_1 d\tau_2 \dots d\tau_n.  \nonumber
\end{equation}
This follows from standard rules of Laplace transforms (see for example
\cite{Roberts:1966}, page 169).  This can be understood by writing
${\cal{L}}^{-1} = {\cal{L}}^{-1}_{\lambda_n}
{\cal{L}}^{-1}_{\lambda_{n-1}} \dots
{\cal{L}}^{-1}_{\lambda_1}$ and then
computing ${\cal{L}}^{-1}_{\lambda_1}[1/\phi(\lambda_1+\dots + \lambda_n)]$ by viewing
$\lambda_2+\dots+\lambda_n$ as a shift.  Using the fact that $t_1 \leq
t_2 \dots \leq t_n$, the above reduces to 
\begin{eqnarray}
{\cal{L}}^{-1} \left[  \frac{1}{\phi(\lambda_1 + \lambda_2 + \dots +
  \lambda_n)} \widetilde{f}(\lambda_1,\lambda_2,\dots,\lambda_n)
\right] (t_1,\dots,t_n) &=& \int_0^{t_1} f(t_1-\tau_1,t_2-\tau_1,\dots,t_n-\tau_1) U'(\tau_1)d\tau_1 \nonumber \\
&=&\int_0^{t_1} f(t_1-\tau,t_2-\tau,\dots,t_n-\tau) dU(\tau).  \label{e:invert1overphi}
\end{eqnarray}
Using $f(\{\lambda_i\}) = U(\{\lambda_i;m_i\})$ from Theorem \ref{t:inversemoments}, equation
(\ref{e:invert1overphi}) lets one compute any order moments by
successive convolutions.  As an application, we give a general expression for the covariance of
an inverse subordinator in terms of $U$ and $dU$.  

\begin{cor}\label{c:inversecov}
Let $D$ be a L\'{e}vy subordinator with L\'{e}vy exponent $\phi$, and let $E$
be the inverse subordinator of $D$.  For $s,t \geq 0$, the covariance $r$, of $E$ is given by
\begin{equation}\label{e:inversecov}
\mbox{\rm Cov}(E(s),E(t)) = \int_0^{s \wedge t}  \left( U(s -
  \tau) + U(t - \tau) \right) dU(\tau) - U(s) U(t).
\end{equation}

\end{cor}

\begin{pf}
 We have
\begin{eqnarray}
\mbox{\rm Cov}(E(s),E(t)) &=& \mathbb{E}E(s)E(t) -
\mathbb{E}E(s)\mathbb{E}E(t) \\
&=& U(s,t;1,1) - U(s) U(t), \label{e:covfirststep}
\end{eqnarray}
where we have used the notation (\ref{e:momentnotation}).  From Theorem \ref{t:inversemoments}, the Laplace transform of
$U(s,t;1,1)$ is 
\begin{eqnarray*}
\widetilde{U}(\lambda_1,\lambda_2;1,1) &=& \frac{1}{\phi(\lambda_1 +
  \lambda_2)} \left( \widetilde{U}(\lambda_1,\lambda_2;1,0) +
  \widetilde{U}(\lambda_1,\lambda_2;0,1) \right) \\
&=& \frac{1}{\phi(\lambda_1 +
  \lambda_2)} \left( \frac{\widetilde{U}(\lambda_1)}{\lambda_2} +
  \frac{\widetilde{U}(\lambda_2)}{\lambda_1} \right).
\end{eqnarray*}
Now, (\ref{e:invert1overphi}) gives
\begin{equation}
U(s,t;1,1) =\int_0^{s \wedge t} \left(U(s-\tau) + U(t-\tau)
\right) dU(\tau). \nonumber
\end{equation}
This with (\ref{e:covfirststep})  establishes (\ref{e:inversecov}).

\end{pf}

\begin{rem}  

Let $\{D(s), \ s \geq 0\}$ be a strictly increasing
  L\'{e}vy subordinator with inverse $\{E(t), t \geq 0 \}$.  Then
  using Corollary \ref{c:inversecov}, one can show that $E$ has
  stationary and uncorrelated increments if and only if $D(s) = C s$
  for some constant $C>0$.

\end{rem}

\section{Examples}\label{s:exbrief}

Examples are discussed in \cite{veillette1:2008}.  We focus there on the
following three important families of subordinators:

\begin{itemize}

\item Poisson and Compound Poisson processes.  Their L\'{e}vy exponent
  is given by 
\begin{equation*}
\phi(\lambda) = \mu \lambda + \int_0^\infty
  (1-e^{-\lambda x}) \nu(dx),
\end{equation*}
 where $\nu$ is a finite measure on $[0,\infty)$.  

\item ``Mixed'' $\alpha$-stable processes.  Their  L\'{e}vy exponent is given by
 \begin{equation*} 
\phi(\lambda) = \int_0^1 \lambda^\beta dp(\beta).
\end{equation*}
  Here $p$ is
  some probability measure on $(0,1)$.  Notice the $\alpha$-stable
  subordinator corresponds to the choice $p(\beta) = \delta(\alpha - \beta)$.

\item Generalized Inverse Gaussian L\'{e}vy processes, which are
  L\'{e}vy processes whose time 1 distribution is given by the
  Generalized inverse Gaussian distribution.  This family of
  distributions was shown to be infinitely divisible in \cite{barndorffn:1977} and its L\'{e}vy-Khintchine
  representation is derived in \cite{eberlein:2002}.   This family includes the
  gamma process, inverse Gaussian process and the reciprocal gamma process.

\end{itemize}

In some cases, an analytic expression for the mean first-hitting time
$U(t)= \mathbb{E}E(t)$ can be given for the above processes, however, this is usually not
the case.  In \cite{veillette1:2008}, we give numerical methods for
inverting the Laplace transform $\widetilde{U}(\lambda) = (\lambda
\phi(\lambda))^{-1}$.  We test these methods in cases where $U$ can be
computed explicitly.  We then study the three examples above
in detail, focusing in each case on the asymptotic behavior of $U$ and on
computing $U$ and higher order moments numerically.

\appendix{

\section{Appendix: Path Properties of Inverse Subordinators}\label{s:introsection}

We describe here the path of the inverse subordinator.  Figure
\ref{f:exampleplot} illustrates the relationship between a subordinator, $\{D(s)$, $s \geq 0\}$, and its
inverse, $\{E(t)$, $t\geq 0\}$.  Observe that both $\{D(s)$, $s \geq
0\}$ and $\{E(t)$, $t\geq 0\}$ are right continuous.  

\begin{figure}[ht]
\centering
\includegraphics[width=0.8\textwidth]{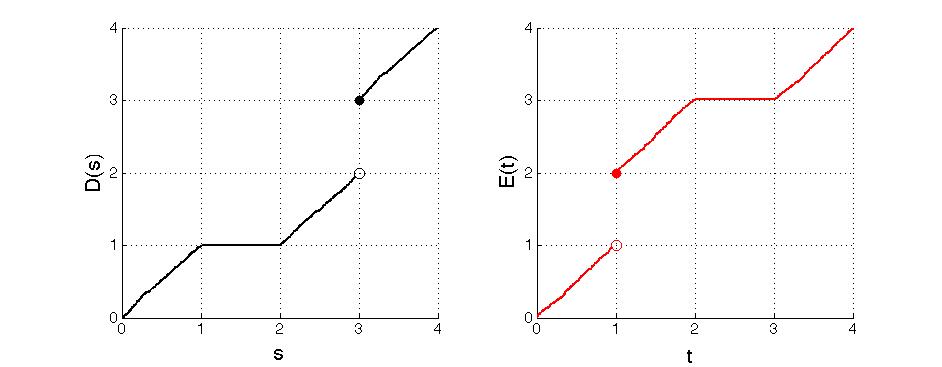}
\caption{An example of a sample path by a subordinator $D$ together
  with its inverse $E$.}
\label{f:exampleplot}
\end{figure}

The following proposition provides additional details.  We provide a proof
for the convenience of the reader.   

\begin{prop}\label{p:eincreasing}
The sample paths of the inverse subordinator $E$ are non-decreasing
and are right continuous with left limits. The sample
paths of $E$ moreover, are continuous if and only if $D$ is strictly increasing\footnote{Sometimes, the inverse is defined as $\inf \{s: D(s) \geq
0 \}$.  This process will have all of the same properties as the
inverse defined here, with the exception that it will be left
continuous as opposed to right continuous.}.
\end{prop}

\begin{pf}

We start by proving that $D$ increasing implies that $E$ is non-decreasing
and has cadlag sample paths (right continuous with left limits).  Indeed, if $t_1 < t_2$,
then $\{s: D(s) > t_2 \} \subset \{s: D(s) > t_1 \}$, meaning $E(t_1)
\leq E(t_2)$.  To prove the sample paths are cadlag, note that the
left limits follow from the fact that $E$ is non-decreasing.  To see right continuity, observe that 
\begin{equation}
\{s: D(s) > t \}  = \bigcup_{\tau > t} \{ s: D(s) > \tau \}, \nonumber
\end{equation} 
since $D$ is non-decreasing.  Thus, taking inf's,
\begin{eqnarray}
E(t)  &=& \inf_s \{ s: D(s) > t \} \nonumber \\
&=&  \inf_s  \bigcup_{\tau > t} \{ s: D(s) > \tau \} \nonumber \\
& = &  \inf_{\tau>t} \inf_s  \{ s: D(s) > \tau \} \quad \mbox{since }
\inf_s  \{ s: D(s) > \tau \}  \in \mathbb{R} \ \mbox{for each} \ \tau, \nonumber \\
& = & \lim_{\tau \rightarrow t^+} E(\tau), \quad \mbox{since}\ \inf_s
\{ s: D(s) > \tau \} \  \mbox{increases in $\tau$}. \nonumber
\end{eqnarray}
 Thus, $\lim_{\tau \rightarrow t^+} E(\tau) =  E(t)$, implying right continuity. 

To prove the second part of the proposition, we will first show that if $E$ is
continuous, then $D$ is strictly increasing.  Suppose the opposite is
true, that is, $D$ is not strictly increasing.  If $D(t) = c$ for $t_1 \leq t
\leq t_2$ with $t_1<t_2$, then by the right continuity of $D$, $E(c-) \leq t_1 $ and $E(c+) > t_2$,
contradicting the continuity of $E$. 

We now show that $D$ strictly increasing implies that $E$ is
continuous.  Suppose that this is not that case, that $E$ is
discontinuous at some point $c$.  From above, $E$ is continuous from the
right, thus discontinuity implies that $E$ is not left-continuous at $c$, meaning $E(c) - E(c-) >
0$.  Then, for any $t$ with $E(c-) < t < E(c) = \inf \{s:D(s) > c\}$,
one has $D(t) \leq c$.  Since $E$ is
non-decreasing, for
$c' < c$ one has 
\begin{equation}
t \geq E(c-) \geq E(c') = \inf \{ s: D(s) > c' \}, \nonumber
\end{equation}
 thus $D(t) >
c'$.  Letting $c' \rightarrow c$ proves that $D(t) =c $, contradicting
that $D$ is strictly increasing.  This finishes the proof.          

 \end{pf}

The next proposition displays the inverse relationship between the
L\'{e}vy subordinator and its first passage time.

\begin{prop}\label{p:setequal}

 Fix $s_1,s_2,\dots, s_n \geq 0$. Then 
\begin{equation}\label{e:probequal}
P[D(s_i) < t_i,\ i=1,\dots,n] = P[E(t_i) > s_i, \ i=1,\dots,n], \quad
\mbox{for a.e.} \ t_1,\dots,t_n \geq 0.
\end{equation}
Moreover, if $D$ is strictly increasing, then (\ref{e:probequal}) holds for all
$t_1,\dots,t_n \geq 0$.

\end{prop}
\begin{pf}
 
First, observe that we have the following set inclusions:
\begin{equation}\label{e:setinclusion}
\{D(s_i) < t_i,\ i=1,\dots,n \} \subset \{ E(t_i) > s_i,\ i=1,\dots,n \} \subset \{D(s_i) \leq t_i,\ i=1,\dots,n \}. 
\end{equation}
To see this, suppose that $D(s_i) < t_i$.  Then by right continuity,
$D(s) < t_i$ for $s > s_i$ sufficiently close to $s_i$.  Thus, $E(t_i)
> s_i$ since $D$ is non-decreasing.  

For the second inclusion, assume $E(t_i) > s_i$.  If $D(s_i) > t_i$,
then we would have $E(t_i) = \inf\{s:D(s) > t_i \} \leq s_i$, a
contradiction.  This verifies (\ref{e:setinclusion}).

Thus, we have for all  $s_1,s_2,\dots, s_n \geq 0$ and  $t_1,\dots,t_n
\geq 0$, we have
\begin{equation}\label{e:inequality}
P[D(s_i) < t_i,\ i=1,\dots,n ] \leq P[ E(t_i) > s_i,\ i=1,\dots,n ]
\leq  P[D(s_i) \leq t_i,\ i=1,\dots,n ]. 
\end{equation}

Now, for  $s_1,s_2,\dots, s_n \geq 0$ fixed,  
\begin{equation}\label{e:otherinequality}
 P[D(s_i) \leq t_i,\ i=1,\dots,n ] \leq P[D(s_i) < t_i,\ i=1,\dots,n ] +
 \sum_{i=1}^n P[D(s_i) = t_i]. 
\end{equation}
For each $i$, let $A_i = \{t: P[D(s_i) = t]  >0 \}$.  Notice that each
$A_i$ is at most a countable subset of $\mathbb{R}^+$, and hence $A
 \equiv \{ (t_1,\dots,t_n) : \ t_i \in A_i \ \mbox{for some} \ i \}$ is a set of Lebesgue
measure $0$.  Thus, (\ref{e:otherinequality}) implies that 
\begin{equation}
 P[D(s_i) \leq t_i,\ i=1,\dots,n ] \leq P[D(s_i) < t_i,\ i=1,\dots,n
 ], \ \ \mbox{for} \  (t_1,\dots,t_n)    \in A^c \nonumber
\end{equation}
This combined with (\ref{e:inequality}) establishes
(\ref{e:probequal}) for all $(t_1,\dots,t_n) \in A^c$.

To prove the second statement of the proposition, assume $D$ is strictly
increasing.  With this, the first inclusion in (\ref{e:setinclusion})
is strengthened to equality:
\begin{equation}\label{e:otherway}
\{D(s_i) < t_i,\ i=1,\dots,n \} = \{ E(t_i) > s_i,\ i=1,\dots,n \}.
\end{equation}
To see this, we proceed by contra-positive.  If $D(s_i) \geq t_i$, then
since $D$ is strictly increasing, $D(s) > t_i$ for all $s > s_i$,
hence $E(t_i) = \inf \{ s: D(s) > t_i \} \leq s_i$. This, combined
with (\ref{e:setinclusion}) proves
(\ref{e:otherway}), thus $P[D(s_i) < t_i,\ i=1,\dots,n ] = P[E(t_i) >
s_i,\ i=1,\dots,n ]$ for all $t_1,\dots,t_n \geq 0$.



\end{pf}

}

\noindent Mark Veillette \& Murad Taqqu \\
\noindent \small Dept. of Mathematics \\
\noindent \small Boston University \\
\noindent \small 111 Cummington St.
\noindent \small Boston, MA 02215

\end{document}